Numerische
Mathematik

# Efficient algorithms for solving the *p*-Laplacian in polynomial time

**Sébastien Loisel[1]**

## Abstract

The *p*-Laplacian is a nonlinear partial differential equation, parametrized by $p \in [1, \infty]$. We provide new numerical algorithms, based on the barrier method, for solving the *p*-Laplacian numerically in $O(\sqrt{n} \log n)$ Newton iterations for all $p \in [1, \infty]$, where $n$ is the number of grid points. We confirm our estimates with numerical experiments.

## 1 Introduction

Let $\Omega \subset \mathbb{R}^d$. For $1 \leq p < \infty$, the $p-$Laplace equation is

$$\nabla \cdot (\|\nabla v\|_2^{p-2} \nabla v) = f \text{ in } \Omega \text{ and } v = g \text{ on } \partial\Omega, \tag{1}$$

where $\|w\|_2 = \left( \sum_{j=1}^d |w_j|^2 \right)^{1/2}$ is the usual $2-$norm on $\mathbb{R}^d$. Prolonging $g$ from $\partial\Omega$ to the interior $\Omega$ and setting $u = v - g$, the variational form is

$$\text{Find } u \in W_0^{1,p}(\Omega) \text{ such that } J(u) = \frac{1}{p} \int_\Omega \|\nabla(u + g)\|_2^p - \int_\Omega fu \text{ is minimized.} \tag{2}$$

A similar definition can be made in the case $p = \infty$ and will be discussed in Sect. 3.1.

For $p = 1$, the *p*-Laplacian is also known as Mean Curvature, and a solution with $f = 0$ is known as a minimal surface [31]. The 1-Laplacian is related to a certain "pusher-chooser" game [19] and compressed sensing [7]. The general *p*-Laplacian is

✉ Sébastien Loisel
  S.Loisel@hw.ac.uk

1 Department of Mathematics, Heriot-Watt University, Riccarton EH14 4AS, UK

Springer

used for nonlinear Darcy flow [11], modelling sandpiles [2] and image processing [8]. We also mention the standard text of Heinonen et al. [16]; as well as the lecture notes of Lindqvist [21].

One may discretize the variational form (2) using finite elements; we briefly outline this procedure in Sect. 2.1 and refer to Barrett and Liu [3] for details. One chooses piecewise linear basis functions $\{\phi_j(x)\}$ on $\Omega$ and we let $u_h(x) = \sum_j u_j \phi_j(x)$. The energy $J(u_h)$ can be approximated by quadrature; the quadrature is exact if the elements are piecewise linear. This leads to a finite-dimensional energy functional

Find $u \in \mathbb{R}^n$ such that $J(u) = c^T u$

$$+ \frac{1}{p} \sum_{i=1}^{m} \omega_i \left( \sum_{j=1}^{d} (D^{(j)}u + b^{(j)})_i^2 \right)^{\frac{p}{2}} \text{ is minimized,} \qquad (3)$$

where $D^{(j)}$ is a numerical partial derivative, $b^{(j)} = D^{(j)}g$ comes from the boundary conditions $g$ and $c$ comes from the forcing term $f$.

Several algorithms have been proposed to minimize the convex functional $J(u)$. Huang et al. [18] proposed a steepest descent algorithm on a regularized functional $J_{h,\epsilon}(u)$ which works well when $p > 2$. Tai and Xu [36] proposed a subspace correction algorithm which works best when $p$ is close to 2 but whose convergence deteriorates when $p \to 1$ or $p \to \infty$. Algorithms based on a multigrid approach (e.g. Huang et al. [9]) suffer from the same problems when $p$ approaches 1 or $\infty$. The algorithm of Oberman [30] also works for $p \geq 2$, although the convergence factor deteriorates after several iterations so it is difficult to reach high accuracy with this method.

The problem of minimizing $J(u)$ has much in common with the problem of minimizing a $p$-norm, which is by now well-understood. The motivation for optimizing a $p$-norm is often given as a facility location problem [1,6]. Efficient algorithms for solving such problems can be obtained within the framework of convex optimization and barrier methods; see Hertog et al. [17] and Xue and Ye [38] specifically for $p$-norm optimization; and for general convex optimization, see Nesterov and Nemirovskii [28], Boyd and Vandenberghe [5] and Nesterov [27] and references therein.

Given a $\nu$-self-concordant barrier for a convex problem, it is well-known that the solution can be found in $O(\sqrt{\nu} \log \nu)$ Newton iterations. However, the "hidden constant" in the big-O notation depends on problem parameters, including the number $n$ of grid points in a finite element discretization. Our main result is to estimate these hidden constants and show that the overall performance of our algorithm is indeed $O(\sqrt{n} \log n)$.

**Theorem 1** *Assume that $\Omega \subset \mathbb{R}^d$ is a polytope and that $T_h$ is a quasi-uniform triangulation of $\Omega$, parametrized by $0 < h < 1$ and with quasi-uniformity parameter $1 \le \rho < \infty$. Let $1 \le p < \infty$. Assume $g \in W^{1,p}(\Omega)$ is piecewise linear on $T_h$ and let $V_h \subset W_0^{1,p}(\Omega)$ be the piecewise linear finite element space on $T_h$ whose trace vanishes. Let $R \ge R^* := 2(1 + \|g\|_{X^p}^p)$, where $|\Omega|$ is the volume of $\Omega$ and*

$$\|g\|_{X^p}^p = \int_\Omega \|\nabla g\|_2^p \, dx.$$

*Let $\epsilon > 0$ be a tolerance. In exact arithmetic, the barrier method of Sect. 2.3, with barrier (25), to minimize $J(u)$ over $u \in V_h$, starting from $(u, s) = (0, \hat{s})$ given by (31), converges in at most $N^*$ iterations, where*

$$N^* \le 14.4\sqrt{|\Omega| h^{-d} d!} \left[ \log\left(h^{-1-17d} R^5 (1 + \|g\|_{X^p}^p) \epsilon^{-1}\right) + K^* \right]. \tag{4}$$

*Here, the constant $K^* = K^*(\Omega, \rho)$ depends on the domain $\Omega$ and on the quasi-uniformity parameter $\rho$ of $T_h$. At convergence, $u$ satisfies*

$$J(u) \le \min_{\substack{v \in V_h \\ \frac{1}{p}\|v + g\|_{X^p}^p \le R}} J(v) + \epsilon.$$

*The global minimizer in $V_h$ can be found by choosing a sufficiently large value of $R$. We give two cases where $R$ is sufficiently large ($1 < p < \infty$ and $p = 1$.)*
**Case $1 < p < \infty$:** *For any $1 < p < \infty$, assume that $f \in L^q(\Omega)$ where $\frac{1}{p} + \frac{1}{q} = 1$. Assume that $\Omega$ fits inside a strip of width $L$. The value of $R = R_{1<p<\infty} = 2 + 8\|g\|_{X^p}^p + 4L^q \left(\frac{p}{2}\right)^{\frac{1}{1-p}} (p-1)\|f\|_{L^q}^q$ always produces the minimizer $u$ of the energy $J(u)$ in the finite element space, and the number of iterations is bounded by*

$$N_{1<p<\infty} \tag{5}$$
$$\le 14.4\sqrt{|\Omega| h^{-d} d!} \left[ \log\left(h^{-1-17d}(1 + \|g\|_{X^p}^p + pL^q\|f\|_{L^q}^q)^5 \epsilon^{-1}\right) + K^* \right].$$

**Case $p = 1$:** *Assume that $L\|f\|_{L^\infty} < 1$. The choice $R = R_{p=1} = 2 + 2\|g\|_{X^1}/(1 - L\|f\|_{L^\infty})$ always produces a global minimizer $u$ of $J(u)$ and the number of iterations is bounded by*

$$N_{p=1} \le 14.4\sqrt{|\Omega| h^{-d} d!} \left[ \log\left(h^{-1-17d}\left(2 + \frac{2\|g\|_{X^1}}{1 - L\|f\|_{L^\infty}}\right)^5 \epsilon^{-1}\right) + K^* \right]. \tag{6}$$

**Computational complexity:** *When we vary $0 < h < 1$ while freezing all other parameters, the three estimates (4), (5) and (6) are $O(\sqrt{n} \log n)$, where $n$ is the number of grid points in $T_h$.*

The iteration count $O(\sqrt{n} \log n)$ also holds if $\epsilon$ is not frozen, provided that $\epsilon^{-1}$ grows at most polynomially in $n$.

We emphasize that the $p = 1, \infty$ cases have up to now been considered to be especially hard and there are no other algorithm that offers any performance guarantees in these situations. Estimate (6) is the first algorithm that is known to converge even when $p = 1$. We also have an algorithm for the $p = \infty$ case and we provide a corresponding iteration count estimate in Theorem 3.

The algorithm mentioned in Theorem 1 is the barrier method of convex optimization. Consider the problem of minimizing $c^T x$ subject to $x \in Q$ where $Q$ is some convex set. Assume we have a $\nu$-self-concordant barrier $F$ for $Q$. The barrier method works by minimizing $t c^T x + F(x)$ for larger and larger values of the barrier parameter $t$, which is increased iteratively according to some schedule. In the short step variant, $t$ increases slowly and the method is very robust; the estimates of Theorem 1 are for the short step variant of the barrier method. It is well-known that the short step barrier method has the theoretically best convergence estimates, but that long step variants (where $t$ increases more rapidly) work better in practice. However, long step algorithms have theoretically worse convergence estimates and, as we will see in the numerical experiments, can sometimes require a large number of Newton iterations to converge.

In order to get the best convergence, we have devised a new, very simple adaptive stepping algorithm for the barrier method. There are already many adaptive stepping algorithms (see e.g. Nocedal et al. [29] and references therein). It is often difficult to prove "global convergence" of these algorithms, and we are not aware of global estimates of Newton iteration counts. With our new, highly innovative algorithm, we are able to prove "quasi-optimal" convergence of our adaptive scheme. Here, quasi-optimal means that our adaptive algorithm requires $\tilde{O}(\sqrt{n})$ Newton iterations, neglecting logarithmic terms, which is the same as the theoretically optimal short-step algorithm of Theorem 1.

The $p$-Laplacian is subject to roundoff problems when $p$ is large but finite, as we now briefly describe. Consider the problem of minimizing $\|v\|_p^p$ in the space $\{v = (1, y)\}$. Assume that we are given a machine $\epsilon$ (for example, $\epsilon \approx 2.22 \times 10^{-16}$ in double precision) and consider an arbitrary vector $v = [1, \delta]^T$. In this situation $\|v\|_p^p = 1^p + \delta^p = 1$ in machine precision, provided $\delta < \epsilon^{1/p}$. This means that a region of size $\epsilon^{1/p}$ near the minimum is numerically indistinguishable from the true minimum when computing the energy, causing a very large relative error in the solution. This phenomenon becomes worse in higher dimensions and when composing with matrices with large condition numbers as in (3). This means that all algorithms, including our own, will struggle to produce highly accurate solutions when $p \gg 2$. In particular, for $p = 5$, we see that $\delta < 7.4 \times 10^{-4}$ is best possible, and this is made worse by the condition number of the differential matrices. However, although the problem is numerically challenging for finite $p \gg 2$, the problem becomes easy again when $p = \infty$. Our second main result is an estimate for the $p = \infty$ case in Sect. 3.1, and we confirm by numerical experiments that there are no numerical issues for $p = \infty$.

Our algorithm is an iterative scheme for a high-dimensional problem arising from a partial differential equation. Each iteration involves the solution of a linear problem that can be interpreted as a numerical elliptic boundary value problem. One can estimate pessimistically that solving each linear problem requires $O(n^3)$ FLOPS, for a total cost of $O(n^{3.5} \log n)$ FLOPS for our entire algorithms. This estimate can be improved by using an $O(n^{2.373})$ FLOPS fast matrix inverse [20], making our overall algorithms $O(n^{2.873} \log n)$ FLOPS; we mention that this matrix inversion algorithm mostly of theoretical interest since it is not practical for any reasonable value of $n$. We have taken special care to preserve the sparsity of this problem so that, if one assumes a bandwidth of $b$ (e.g. typically $b = O(\sqrt{n})$ for $d = 2$ and $b = O(n^{2/3})$ for $d = 3$), one obtains an $O(b^2 n)$ sparse matrix solve algorithm, resulting in $O(n^{2.5} \log n)$ ($d = 2$) or $O(n^{2.84} \log n)$ ($d = 3$) FLOPS for our overall algorithms. In addition, we mention many preconditioning opportunities [4,10,12,14,15,22–26,35]. Although solution by preconditioning is possible, it is difficult to estimate the number of iterations a priori since the diffusion coefficient of the stiffness matrix is difficult to estimate a priori; in the best case ("optimal preconditioning") where the elliptic solve at each Newton iteration can be done in $O(n)$ FLOPS, our algorithms are then $O(n^{1.5} \log n)$ FLOPS.

Our paper is organized as follows. In Sect. 2, we give some preparatory material on the *p*-Laplacian and the barrier method. In Sect. 3, we prove our main theorem for $1 \leq p < \infty$ and a separate theorem for the case $p = \infty$. In Sect. 4, we validate our algorithms with numerical experiments. We end with some conclusions.

## 2 Preparatory material

We now discuss some preparatory material regarding the $p-$Laplacian, for $1 \leq p \leq \infty$. The $\infty$-Laplacian can be interpreted as the problem of minimizing

$$J(u) = \|u + g\|_{X^\infty(\Omega)} - \int_\Omega fu \text{ where } \|v\|_{X^\infty(\Omega)} = \sup_{x \in \Omega} \|\nabla v(x)\|_2. \qquad (7)$$

Note that (7) is not a limit as $p \to \infty$ of (2), e.g. because (2) uses the $p$th power of $\| \cdot \|_{X^p}$ in its definition.

**Lemma 1** *For $1 \leq p \leq \infty$, $J(u)$ is convex on $W^{1,p}(\Omega)$. For $1 < p < \infty$, $J(u)$ is strictly convex on $W_0^{1,p}(\Omega)$.*

***Proof*** We consider the case $1 \leq p < \infty$ in detail, the case $p = \infty$ is done in a similar fashion. Convexity (and strict convexity) is unaffected by linear shifts so without loss of generality we assume that $f = 0$. Let $0 \leq t \leq 1$. We must show that $J(tu+(1-t)v) \leq t J(u)+(1-t) J(v)$. To simplify the notation, let $q = \nabla u, r = \nabla v$ and $s = \nabla g$.

$$J(tu + (1 - t)v) = \int_\Omega \|tq + (1 - t)r + s\|_2^p \overset{(*)}{\leq} \int_\Omega (t\|q + s\|_2 + (1 - t)\|r + s\|_2)^p$$

$$\overset{(**)}{\leq} \left( \int_\Omega t\|q + s\|_2^p + (1 - t)\|r + s\|_2^p \right) = t J(u) + (1 - t) J(v),$$

where we have used the triangle inequality for $\| \cdot \|_2$ at $(*)$ and the convexity of $\phi(z) = z^p$ at $(**)$.

We now prove strict convexity for the $1 < p < \infty$ case. If we have equality at $(*)$ then $q(x) + s(x)$ and $r(x) + s(x)$ are non-negative multiples of one another, i.e. $q + s = aw$ and $r + s = bw$ where $a(x), b(x) \geq 0$ and $w(x)$ is vector-valued. Then $(**)$ becomes $\int_\Omega ((ta + (1-t)b)\|w\|_2)^p \leq \int_\Omega (ta^p + (1-t)b^p)\|w\|_2^p$. Note that $(ta+(1-t)b)^p < ta^p + (1-t)b^p$ unless $a = b$ so the inequality $(**)$ is strict unless $\nabla u = \nabla v$ almost everywhere. Since $u, v \in W_0^{1,p}$ can be identified by their gradients, we have proven strict convexity. □

From the norm equivalence $\|u\|_p \leq d^{\max\left\{0, \frac{1}{p} - \frac{1}{q}\right\}} \|u\|_q$ for $x \in \mathbb{R}^d$, one obtains

$$d^{-\max\left\{0, \frac{1}{p} - \frac{1}{2}\right\}} |u|_{W^{1,p}} \leq \|u\|_{X^p(\Omega)} \leq d^{\max\left\{0, \frac{1}{2} - \frac{1}{p}\right\}} |u|_{W^{1,p}}. \quad (8)$$

We can give a modified Friedrichs inequality for $\|\cdot\|_{X^p}$.

**Lemma 2** (Friedrichs inequality for $\|\cdot\|_{X^p}$) *Assume that $\Omega \subset \mathbb{R}^d$ fits inside of a strip of width $\ell$ and assume that $\phi \in W_0^{1,p}(\Omega)$, where $1 \leq p \leq \infty$. Then, $\|\phi\|_{L^p} \leq \ell p^{-\frac{1}{p}} \|\phi\|_{X^p}$, where we define $\infty^{-\frac{1}{\infty}} = 1$.*

**Proof** Without loss of generality, assume that $\Omega$ is inside the strip $0 \leq x_1 \leq \ell$. From the fundamental theorem of calculus, the following argument proves the $p = \infty$ case: $|\phi(x_1, \ldots, x_d)| \leq \int_0^{x_1} |\phi_{x_1}(t, x_2, \ldots, x_d)| \, dx \leq \int_0^\ell \sup_{x \in \Omega} \|\nabla \phi(x)\|_2 \, dx_1 \leq \ell \|\phi\|_{X^\infty}$. Now assume $1 \leq p < \infty$.

$$\int_0^\ell |\phi|^p \, dx_1 = \int_0^\ell \left| \int_0^{x_1} \phi_{x_1}(t, x_2, \ldots, x_d) \, dt \right|^p dx_1 \quad (9)$$

$$\leq \int_0^\ell \int_0^{x_1} |\phi_{x_1}(t, x_2, \ldots, x_d)|^p x_1^{p-1} \, dt \, dx_1 \text{ (Jensen's ineq.)} \quad (10)$$

$$\leq \int_0^\ell \|\nabla \phi(t, x_2, \ldots, x_d)\|_2^p \int_t^\ell x_1^{p-1} \, dx_1 \, dt \quad (11)$$

$$= \int_0^\ell \|\nabla \phi(t, x_2, \ldots, x_d)\|_2^p \frac{1}{p} (\ell^p - t^p) \, dt \quad (12)$$

$$\leq \frac{\ell^p}{p} \int_0^\ell \|\nabla \phi(t, x_2, \ldots, x_d)\|_2^p \, dt. \quad (13)$$

The result follows by integrating over $x_2, \ldots, x_d$. □

We now give an *a priori* estimate on the magnitude of the minimizer of $J(u)$. This estimate will be important in the design of our algorithm in order to limit the search volume to some ball of reasonable size.

**Lemma 3** *Let $1 < p < \infty$ and $\frac{1}{p} + \frac{1}{q} = 1$ and assume that $\Omega \subset \mathbb{R}^d$ is a domain of width $L$. Let $\|v\|_{X^p}^p = \int_\Omega \|\nabla v\|_2^p$. Assume $\{u_k\} \subset W_0^{1,p}(\Omega)$ is a minimizing sequence for $J(u)$. Then, for large enough $k$,*

$$\|u_k\|_{X^p}^p \leq 4\|g\|_{X^p}^p + 2L^q \left(\frac{p}{2}\right)^{\frac{1}{1-p}} (p-1)\|f\|_{L^q}^q. \quad (14)$$

*If $\{p, q\} = \{1, \infty\}$ and $L\|f\|_{L^q} < 1$ then a minimizing sequence must eventually lie in $\|u_k\|_{X^p} \leq \|g\|_{X^p}/(1 - L\|f\|_{L^q})$.*

**Proof** Case $1 < p < \infty$: For convenience, we write $J(u) = \frac{1}{p}\|u + g\|_{X^p}^p - \int_\Omega fu$. Assume $\|u\|_{X^p} \geq \|g\|_{X^p}$; then:

$$J(u) \geq \frac{1}{p}(\|u\|_{X^p} - \|g\|_{X^p})^p - \|f\|_{L^q}\|u\|_{L^p}$$
$$\geq \frac{1}{p}\|u\|_{X^p}^p - \frac{1}{p}\|g\|_{X^p}^p - \|f\|_{L^q}Lp^{-1/p}\|u\|_{X^p}.$$

Next, we use Young's inequality $ab \leq \frac{1}{q}a^q + \frac{1}{p}b^p$ with $a = 2^{1/p}p^{-1/p}L\|f\|_{L^q}$, $b = 2^{-1/p}\|u\|_{X^p}$ to obtain

$$J(u) - J(0) \geq \frac{1}{2p}\|u\|_{X^p}^p - \frac{2}{p}\|g\|_{X^p}^p - \frac{L^q(p/2)^{\frac{1}{1-p}}}{q}\|f\|_{L^q}^q. \qquad (15)$$

Hence, if $\|u\|_{X^p}^p > 4\|g\|_{X^p}^p + 2L^q(p/2)^{\frac{1}{1-p}}(p-1)\|f\|_{L^q}^q$, then $J(u) - J(0) > 0$ and hence a minimizing sequence must satisfy (14).

The $p = 1$ case is as follows:

$$J(u) - J(0) \geq \|u\|_{X^1} - \|g\|_{X^1} - \|f\|_{L^\infty}L\|u\|_{X^1} > 0,$$

if $\|u\|_{X^1} > \|g\|_{X^1}/(1 - \|f\|_{L^\infty}L)$. The $p = \infty$ case is done in a similar fashion. $\quad\square$

The *a priori* estimate above can also be used to show the existence of a minimizer of $J(u)$.

**Lemma 4** *Let $1 < p < \infty$. There is a unique $u \in V \subset W_0^{1,p}(\Omega)$ that minimizes $J(u)$.*

**Proof** Let $\alpha = \inf_v J(v)$. We now show how to produce a minimizing sequence $\{u_k\} \subset W_0^{1,p}(\Omega)$. For $k = 1, 2, \ldots$, let $B_k = \{u \in W_0^{1,p}(\Omega) \mid J(u) < \alpha + 1/k$ and $\|u\|_{X^p} < 4\|g\|_{X^p}^p + 4L^2(p-1)\|f\|_{L^q}^q + 1\}$, see (14). Note that each $B_k$ is open and nonempty so pick $u_k \in B_k$. Furthermore, the $B_k$ are nested: $B_1 \supset B_2 \supset \ldots$; the convexity of $J$ implies that the $B_k$ are also convex.

According to (8), we see that each $B_k$ is contained in a closed ball $F = \{u \in W_0^{1,p}(\Omega) \mid |u|_{W^{1,p}(\Omega)} \leq r\}$ where $r = d^{\max\left\{0, \frac{1}{p} - \frac{1}{2}\right\}}\left(4\|g\|_{X^p}^p + 4L^2(p-1)\|f\|_{L^q}^q + 1\right)$. Recall that $F$ is weakly compact. Passing to a subsequence if necessary, we may assume that $\{u_k\}$ converges weakly to some $u$. By Mazur's lemma, we can now find some convex linear combinations $v_k = \sum_{j=k}^{J(k)}\alpha_j u_j \in B_k$ such that $\{v_k\}$ converges to $u$ strongly. This shows that $u$ belongs to every $B_k$ and hence $J(u)$ is minimal. Uniqueness follows by strict convexity. $\quad\square$

## 2.1 Finite elements

Assume that $\Omega$ is a polygon. We introduce a triangulation $T_h$ of $\Omega$, parametrized by $0 < h < 1$, and piecewise linear finite element basis functions $\{\phi_1(x), \ldots, \phi_n(x)\} \subset W^{1,p}(\Omega)$. As usual, we define a "reference element" $\hat{K} = \{x \in \mathbb{R}^d : x_i \geq 0 \text{ for } i = 1, \ldots, d \text{ and } \sum_{i=1}^{d} x_i \leq 1\}$. Each simplex $K_k$ in $T_h$ can be written as $K_k = \Phi_k(\hat{K}) = P^{(k)}\hat{K} + q^{(k)}$, where $P^{(k)} \in \mathbb{R}^{d \times d}$ and $q^{(k)} \in \mathbb{R}^d$. If $T_h$ is a uniform lattice of squares or $d$-cubes, then each $P^{(k)}$ is of the form $\text{diag}(\pm h, \ldots, \pm h)$, and $\|\|P^{(k)}\|\|_2 = \|\|[P^{(k)}]^{-1}\|\|_2^{-1} = h$. In general, if $T_h$ is not necessarily a uniform lattice, we say that the family of triangulations $T_h$, parametrized by $0 < h < 1$, is **quasi-uniform** with parameter $\rho < \infty$ if $h \leq \left\{\|\|P^{(k)}\|\|_2, \|\|[P^{(k)}]^{-1}\|\|_2^{-1}\right\} \leq \rho h$. Note that on the reference simplex, the basis functions are $\hat{\phi}_i(\hat{x}) = \hat{x}_i$ for $i = 1, \ldots, d$ and $\hat{\phi}_0(\hat{x}) = 1 - \sum_i \hat{x}_i$. As a result, $\|\nabla\hat{\phi}\|_2 \leq \sqrt{d}$ and, from the chain rule, $\|\nabla\phi_i(x)\|_2 \leq h^{-1}\sqrt{d}$.

Let $\text{span}\{\phi_k(x) \mid k = 1, \ldots, n\} \subset W^{1,p}(\Omega)$ be the finite element space of piecewise linear elements over $T_h$ and let $\int_\Omega w(x)\,dx \approx \sum_{i=1}^{m} \omega_i w(x^{(i)})$ be the midpoint quadrature rule, which is exact for piecewise linear or piecewise constant functions. We can construct a "discrete derivative" matrix $D^{(j)}$ whose $(i, k)$ entry is $D_{i,k}^{(j)} = \frac{\partial\phi_k}{\partial x_j}(x^{(i)})$. Then,

$$\frac{1}{p}\int_\Omega \|\nabla(u + g)\|_2^p = \sum_{i=1}^{m} \frac{\omega_i}{p}\left(\sum_{j=1}^{d}\left((D^{(j)}(u + g))_i\right)^2\right)^{\frac{p}{2}};$$

note that the quadrature is exact provided that $g$ is also piecewise linear. For the midpoint rule, $\omega_i$ is the volume of the simplex $K_i$; if the triangulation $T_h$ is quasi-uniform then we find that

$$\frac{h^d}{d!} \leq \omega_i \leq \frac{\rho^d h^d}{d!}; \tag{16}$$

we write $\omega_i = \Theta(h^d)$, which means both that $\omega_i = O(h^d)$ and $h^d = O(\omega_i)$. We abuse the notation and use the same symbol $u$ to represent both the finite element coefficient vector $[u_1, \ldots, u_n]^T$ and the finite element function $u(x) = \sum_{k=1}^{n} u_k\phi_k(x)$.

We further denote by $D_\Gamma^{(j)}$ the columns of $D^{(j)}$ corresponding to the vertices of $T_h$ in $\partial\Omega$, and $D_I^{(j)}$ corresponding to the interior vertices in $\Omega$, such that $D^{(j)} = \left[D_I^{(j)} \; D_\Gamma^{(j)}\right]$. Denoting $u = \left[u_I^T \; u_\Gamma^T\right]^T = \left[u_I^T \; 0^T\right]^T$, note that $D^{(j)}(u + g) = D_I^{(j)}u_I + D^{(j)}g$. Putting $b^{(j)} = D^{(j)}g$ and dropping the subscript $I$ leads to the discretized system (3). The matrix $A = \sum_{k=1}^{d}[D_I^{(k)}]^T \text{diag}(\omega_1, \ldots, \omega_m)D_I^{(k)}$ is the usual discretization of the Laplacian or Poisson problem, and we have that $u^T A u = |u|_{H^1}^2 = \int_\Omega \|\nabla u\|_2^2\,dx$. For a domain of width $\ell$, the Friedrichs inequality $\|u\|_{L^2} \leq \ell|u|_{H^1}$ (see [33, (18.1) and (18.19)]) proves that the smallest eigenvalue of the Laplacian differential operator is at least $\ell^{-2}$; however, the smallest eigenvalue of the finite-dimensional matrix $A$ is actually $\Theta(h^2)$ because the relevant Rayleigh quotient in $\mathbb{R}^n$ is $u^*Au/u^*u \neq |u|_{H^1}^2/\|u\|_{L^2}^2$.

We now prove that the finite element method converges for the *p*-Laplacian.

**Lemma 5** *Assume that $\Omega$ is a polytope and $1 < p < \infty$. Let $u_h$ be the finite element minimizer of $J(u)$ in a finite element space that contains the piecewise linear functions. Then, $J(u_h)$ converges to $\inf_v J(v)$ as $h \to 0$.*

**Proof** Let $u$ be a minimizer of $J(u)$ and denote by $V_h \subset W_0^{1,p}(\Omega)$ the finite element space with grid parameter $0 < h < 1$. Recall that finite element functions are dense in $W_0^{1,p}(\Omega)$ [13, Proposition 2.8, page 316]. Hence, we can find finite element functions $\{v_h \in V_h\}$ that converge to $u$ in the $W_0^{1,p}(\Omega)$ norm as $h \to 0$. Since $J$ is continuous and since $u_h$ minimizes $J(u)$ in the finite element space $V_h$, we find that $J(u_h) \leq J(v_h) \to J(u) = \inf_v J(v)$, as required. □

Lemma 5 is very general (no regularity assumptions are made on the solution $u$) but also very weak since it does not give a rate of convergence. If one assumes regularity of the solution then one can use quasi-interpolation [34] to estimate the convergence more precisely. However, we will see in Sect. 2.2 (Example 1) that it is difficult to prove regularity. Since the present paper focuses on the numerical solver, and not in the discretization, we do not investigate this aspect any further. The theorem also does not specify whether $u_h$ converges as $h$ tends to 0. In the case $1 < p < \infty$, the strict convexity of $J$ ensures that $u_h$ will indeed converge to a $u$ in $W^{1,p}(\Omega)$ but for $p = 1$ there may be multiple minimizers and then $u_h$ could oscillate between the many minimizers or converge to a "minimizer" in the double-dual of $W^{1,1}(\Omega)$.

## 2.2 Pathological situations for extreme *p* values

The *p*-Laplace problem varies in character as $p$ ranges from $1 \leq p \leq +\infty$. When $p = 2$, minimizing $J(u)$ is equivalent to solving a single linear problem, which is clearly faster than solving hundreds of linear problems as required by a barrier method. As $p$ gets further away from $p = 2$, naive solvers work less well and proper convex optimization algorithms are required, such as our barrier methods. The extreme cases $p = 1$ and $p = \infty$ have traditionally been considered hardest. For example, $J(u)$ may not be differentiable at $u$ when $p \in \{1, \infty\}$, typically when $\nabla(u+g)$ vanishes at some point $x \in \Omega$. In Sect. 2, we have introduced several lemmas, some of which work for all cases $1 \leq p \leq \infty$, others are restricted to $1 < p < \infty$. Briefly speaking, we have shown that for all $1 \leq p \leq \infty$, $J(u)$ is convex and possesses minimizing sequences (with some restrictions on the forcing $f$ when $p \in \{1, \infty\}$.) These facts are sufficient to deploy barrier methods, because barrier methods do not require the objective to be differentiable, be strictly convex, or have a unique minimizer. As a "bonus", we have also shown that $J(u)$ is strictly convex and has a unique minimum when $1 < p < \infty$, but this is not required for the successful application of our barrier methods.

We now illustrate the pathological behavior for extreme values of $p$ with several simple examples. For $1 < p < \infty$, strict convexity ensures the uniqueness of the minimum of $J(u)$. However, in the case $p = 1$, the minimizer may be "outside" of $W^{1,1}(\Omega)$ or nonunique.

**Example 1** Consider $\Omega = (0, 1)$ in dimension $d = 1$ and $f = 0$, with boundary conditions $u(0) = 0$ and $u(1) = 1$ and with $p = 1$. Then,

$$J(u) = \int_0^1 |u'(x)| \, dx = TV(u) \geq 1,$$

where $TV(u)$ denotes the usual total variation of $u$. Any monotonically nondecreasing function $u(x)$ with $u(0) = 0$ and $u(1) = 1$ will minimize $J(u)$ and satisfy the boundary conditions.

A minimizing sequence for $J(u)$ is the piecewise linear functions $u_n(x) = \min(1, \max(0, 0.5 + n(x - 0.5)))$. This sequence converges to the indicating function of $[0.5, 1)$, which is not in $W^{1,1}(0, 1)$. This is because $W^{1,1}(0, 1)$ is not reflexive and hence its unit ball is not weakly compact. Instead, the limit of $u_n$ is in $BV$, the double-dual of $W^{1,1}(0, 1)$.

We now briefly show why the minimization of $J(u)$ for $u \in V_h$ is numerically challenging.

**Example 2** Consider $J(x) = |x|^p$ where $x \in \mathbb{R}$ and $1 \leq p < \infty$; this corresponds to a 1-dimensional discrete $p$-Laplacian with a single grid point. The Newton iteration $x_{k+1} = x_k - J'(x_k)/J''(x_k)$ is

$$x_{k+1} = x_k - \frac{\text{sgn}(x_k)p|x_k|^{p-1}}{p(p-1)|x_k|^{p-2}} = \frac{p-2}{p-1}x_k.$$

Hence, the Newton iteration converges linearly for $p \in (1.5, 2) \cup (2, \infty)$ and diverges for $1 < p \leq 1.5$. The Newton iteration is undefined for $p = 1$ since $J'' = 0$.

The $p$-Laplacian for $p = 1$ is particularly hard; we now show two types of difficulties. First, the Hessian may be singular, and regularizing the Hessian leads to gradient descent.

**Example 3** Consider $J(x) = \sqrt{x_1^2 + x_2^2} = \|x\|_2$; this correspond to a 2-dimensional 1-Laplacian discretized with a single grid point. The gradient is $J'(x) = \frac{x}{\|x\|_2}$ and the Hessian is

$$J''(x) = \frac{1}{\|x\|_2^3} \begin{bmatrix} x_2^2 & -x_1x_2 \\ -x_1x_2 & x_1^2 \end{bmatrix}.$$

The Hessian matrix $J''(x)$ is singular which makes the Newton iteration undefined. To make matters worse, the kernel of $J''$ is spanned by $J'$ and hence any "regularization" $J'' + \epsilon I$ leads to a simple gradient descent.

Yet another difficulty is that the 1-Laplacian may have nonunique solutions or no solutions when the forcing is nonzero.

**Example 4** Let $c \in \mathbb{R}$ and $J(x) = |x| + cx$; this corresponds to a 1-dimensional 1-Laplacian with a nonzero forcing term, discretized with a single grid point. Then, $J(x)$ is convex for all $c \in \mathbb{R}$. However, $J(x)$ has a unique minimum $x = 0$ if and only if $|c| < 1$. When $|c| = 1$, $J(x)$ has infinitely many minima. When $|c| > 1$, $J(x)$ is unbounded below and there is no minimum.

As a result, the energy $J(u)$ of the 1-Laplacian may not be bounded below when the forcing $f \not\equiv 0$; see also Lemma 3.

## 2.3 Convex optimization by the barrier method

In this section, we briefly review the theory and algorithms of convex optimization and refer to Nesterov [27, Section 4.2] for details, including the notion of self-concordant barriers.

Let $Q \subset \mathbb{R}^n$ be a bounded closed convex set that is the closure of its interior, $c \in \mathbb{R}^n$ be a vector and consider the convex optimization problem

$$c^* = \min\{c^T x \ : \ x \in Q\}. \tag{17}$$

The barrier method (or interior point method) for solving (17) is to minimize $tc^T x + F(x)$ for increasing values of $t \to \infty$, where the barrier function $F(x)$ tends to $\infty$ when $x \to \partial Q$. The minimizer $x^*(t)$, parametrized by $t \geq 0$, of $tc^T x + F(x)$, is called the **central path**, and $x^*(t)$ forms a minimizing sequence[1] for (17) as $t \to \infty$. Assume we have a $\nu$-self-concordant barrier $F(x)$ for $Q$. Define the norm $\|v\|_x^* = \sqrt{v^T [F''(x)]^{-1} v}$. The **main path-following scheme** is

1. Set $t_0 = 0$, $\beta = 1/9$ and $\gamma = 5/36$. Choose an accuracy $\epsilon > 0$ and $x^{(0)} \in Q$ such that $\|F'(x^{(0)})\|_{x^{(0)}}^* \leq \beta$.
2. The $k$th iteration ($k \geq 0$) is

$$t_{k+1} = t_k + \frac{\gamma}{\|c\|_{x^{(k)}}^*} \text{ and } x^{(k+1)} = x^{(k)} - [F''(x^{(k)})]^{-1}(t_{k+1}c + F'(x^{(k)})). \tag{18}$$

3. Stop if $t_k \geq \left(\nu + \frac{(\beta + \sqrt{\nu})\beta}{1-\beta}\right)\epsilon^{-1} =: \text{tol}^{-1}$.

The invariant of this algorithm is that, if $\|t_k c + F'(x^{(k)})\|_{x^{(k)}}^* \leq \beta$ then also $\|t_{k+1}c + F'(x^{(k+1)})\|_{x^{(k+1)}}^* \leq \beta$. The stopping criterion guarantees that, at convergence, $c^T x^{(k)} - c^* \leq \epsilon$. Starting this iteration can be difficult, since it is not always obvious how to find an initial point $x^{(0)} \in Q$ such that $\|F'(x^{(0)})\|_{x^{(0)}}^* \leq \beta$. Define the **analytic center** $x_F^*$ by $F'(x_F^*) = 0$. We use an **auxiliary path-following scheme**[2] to approximate the analytic center $x_F^*$ of $Q$:

1. Choose $x^{(0)} \in Q$ and set $t_0 = 1$ and $G = -F'(x^{(0)})$.
2. For the $k$th iteration ($k \geq 0$):

$$t_{k+1} = t_k - \frac{\gamma}{\|G\|_{x^{(k)}}^*} \text{ and} \tag{19}$$

---

[1] Perhaps one should say "minimizing filter".

[2] Our presentation corrects some misprints in the auxiliary path algorithm of Nesterov [27, Section 4.2].

$$x^{(k+1)} = x^{(k)} - [F''(x^{(k)})]^{-1}(t_{k+1}G + F'(x^{(k)})). \tag{20}$$

3. Stop if $\|F'(x^{(k)})\|_{x^{(k)}}^* \leq \frac{\sqrt{\beta}}{1+\sqrt{\beta}}$. Set $\bar{x} = x^{(k)} - [F''(x^{(k)})]^{-1}F'(x^{(k)})$.

The invariant of the auxiliary scheme is that $\|t_k G + F'(x^{(k)})\|_{x^{(k)}}^* \leq \beta$ for every $k$. At convergence, one can show that $\|F'(\bar{x})\|_{\bar{x}}^* \leq \beta$. Let $\hat{x} \in Q$ be some starting point for the auxiliary path-following scheme. Combining the auxiliary path-following scheme to find the approximate analytic center $\bar{x}$ of $Q$, followed by the main path-following scheme to solve the optimization problem (17) starting from $x^{(0)} = \bar{x}$, completes in at most $N$ iterations, where

$$N = 7.2\sqrt{\nu}\left[2\log\nu + \log\|F'(\hat{x})\|_{x_F^*}^* + \log\|\hat{x}\|_{x_F^*}^* + \log(1/\epsilon)\right]. \tag{21}$$

### 2.3.1 Long-step algorithms

The path-following schemes of Sect. 2.3 are so-called "short step", meaning that the barrier parameter $t$ increases fairly slowly when $\nu$ is large. It is well-known that long-step algorithms, where $t$ increases more rapidly, often converge faster overall than short-step algorithms, even though the worst case estimate $O(\nu\log\nu)$ is worse than the short-step estimate $O(\sqrt{\nu}\log\nu)$, see Nesterov and Nemirovskii [28] for details. The main path-following scheme can be made "long-step" as follows:

1. Assume $x^{(0)} \in Q$ such that $\|F'(x^{(0)})\|_{x^{(0)}}^* \leq \beta$ and let $t_0 = 0$.
2. Set

$$t_{k+1} = \begin{cases} \max\left\{\kappa t_k, t_k + \frac{\gamma}{\|c\|_{x^{(k)}}^*}\right\} & \text{if } \|t_k c + F'(x^{(k)})\|_{x^{(k)}}^* \leq \beta, \\ t_k & \text{otherwise}; \end{cases} \tag{22}$$

$$x^{(k+1)} = x^{(k)} - r_k[F''(x^{(k)})]^{-1}(t_{k+1}c + F'(x^{(k)})), \tag{23}$$

where $0 < r_k \leq 1$ is found by line search, see e.g. Boyd and Vandenberghe [5, Algorithm 9.2 with $\alpha = 0.01$ and $\beta = 0.25$].
3. Stop if $t_k \geq \left(\nu + \frac{(\beta+\sqrt{\nu})\beta}{1-\beta}\right)\epsilon^{-1} = \text{tol}^{-1}$.

The parameter $\kappa \geq 1$ determines the step size of the scheme. In convex optimization, step sizes $\kappa = 10$ or even $\kappa = 100$ are often used, but we will see in Sect. 4 that shorter step sizes are better suited for the $p$-Laplacian.

The long-step variant of the auxiliary path-following scheme is implemented in a similar fashion; the criterion for decreasing $t_{k+1}$ is then $\|t_k G + F'(x^{(k)})\|_{x^{(k)}}^* \leq \beta$.

### 2.3.2 Adaptive stepping

We finally introduce an algorithm whose step parameter $\kappa_k$ is indexed by the iteration counter $k$. We first introduce some terminology. If $\|t_k c + F'(x^{(k)})\|_{x^{(k)}}^* \leq \beta$ (main phase) or $\|t_k G + F'(x^{(k)})\|_{x^{(k)}}^* \leq \beta$ (auxiliary phase), we say that $x^{(k)}$ was **accepted**,

else we say that $x^{(k)}$ was a **slow step**. Let $\kappa_0$ be an initial step size (we will take $\kappa_0 = 10$.)

1. If $x^{(k)}$ is accepted after 2 or fewer slow steps, put $\kappa_{k+1} = \min\{\kappa_0, \kappa_k^2\}$.
2. If $x^{(k)}$ is accepted after 8 or more slow steps, put $\kappa_{k+1} = \sqrt{\kappa_k}$.
3. If $x^{(k)}$ is still not accepted after 15 slow steps, replace $x^{(k+1)}$ and $t_{k+1}$ by the most recently accepted step and put $\kappa_{k+1} = \kappa_k^{1/4}$. We call this procedure a **rejection**.
4. Otherwise, put $\kappa_{k+1} = \kappa_k$.

The quantity $t_{k+1}$ is computed as in the long step algorithm (22), with $\kappa = \kappa_{k+1}$. Note that whenever $t_{k+1}$ coincides with the short step (20), then the step is automatically accepted. The rejection is "wasteful" in that it discards possibly useful information, but we will see in the numerical experiments that this adaptive scheme is quite efficient in practice. Furthermore, the rejection step is the key that unlocks a very simple analysis of our algorithm.

**Theorem 2** *For given $c$, $F$, $\epsilon$, let $N_S$ and $N_A$ be the number of Newton steps of the short step and adaptive step algorithms, respectively. Then,*

$$N_A \le 16\lceil 0.76 + 0.73 \log(1 + 9\sqrt{\nu})\rceil N_S. \tag{24}$$

**Proof** By construction, on each accepted step of the main path-following algorithm, we find that $t_{k+1} \ge t_k + \frac{\gamma}{\|c\|_{x^{(k)}}^*}$, the short step size, see (22). Thus, we only need to estimate the maximum number of slow steps before a step is accepted. According to [27, p.202], the short step size satisfies

$$t_{k+1} \ge \overbrace{\left(1 + \frac{5}{4 + 36\sqrt{\nu}}\right)}^{\kappa_{\min}} t_k.$$

Starting from $\kappa = 10$, after $r$ rejections, the step size is $\kappa = 10^{(1/4)^r}$. When $\kappa \le \kappa_{\min}$, the short step is automatically accepted and hence the maximum number of rejections is $r = \lceil r_-\rceil$, where

$$10^{(1/4)^{r_-}} = \kappa_{\min} \implies r_- = -\frac{\log(\log(\kappa_{\min})/\log(10))}{\log 4}.$$

Hence,

$$r \le \lceil 0.76 + 0.73 \log(1 + 9\sqrt{\nu})\rceil.$$

Since all the adaptive steps are at least as large as the short steps and the stopping criterion is purely based on the barrier parameter $t_k$, and noting that each rejection corresponds to 15 slow steps (plus the initial accepted step), we obtain the estimate for the main phase. The estimate for the auxiliary phase is obtained in a similar fashion. $\square$

Theorem 2 states that the adaptive algorithm cannot be much worse than the short step algorithm, which means that the adaptive algorithm scales at worse like $\tilde{O}(\sqrt{\nu})$, where we have neglected some logarithms. The reader may be surprised that the estimate for the adaptive scheme is slightly worse than the estimate for the short step scheme, but this is a well-known phenomenon in convex optimization. The long step estimates are quite pessimistic and in practice, long step and adaptive schemes work much better than the theoretical estimates. Our result is especially interesting because it is well-known that estimates for long step algorithms scale like $\tilde{O}(\nu)$, whereas our new algorithm scales like $\tilde{O}(\sqrt{\nu})$.

## 3 Proof of Theorem 1

The proof of Theorem 1 is rather technical, so we begin by outlining the plan of our proof. The idea is to estimate all the quantities in the bound (21) for the number $N$ of Newton iterations. The barrier parameter $\nu$ is estimated in Lemma 6. Some "uniform" or "box" bounds are given for the central path in Lemma 7; these are an intermediate step in converting as many estimates as possible into functions of $h$. Because (21) depends on the Hessian $F''$ of the barrier, the lowest eigenvalue $\lambda_{\min}$ of $F''$ is estimated in Lemma 8. This bound itself depends on extremal singular values of the discrete derivative matrices, which are estimated in Lemma 9, and these bounds are rephrased in terms of $h$ in Lemma 10. In Lemma 11, we establish the connection between the number $m$ of simplices and the grid parameter $h$, which is used in Lemma 12 to estimate the quantities $\|\hat{x}\|_2$ and $\|F'(\hat{x})\|_2$, which can be converted to estimates for $\|\hat{x}\|_{x_F^*}^*$ and $\|F'(\hat{x})\|_{x_F^*}^*$ in (21) by dividing by $\lambda_{\min}$; here $\hat{x}$ is a starting point for the barrier method. Finally, the quantities $R$ appearing in Theorem 1 are obtained by starting from the estimates of Lemma 3, adding 1, and doubling them. This ensures that the central path will be well inside the ball of radius $R$.

In the present section, we treat in detail the case $1 \leq p < \infty$. The case $p = 1$, which is considered especially difficult, poses no special difficulty in the present section, provided that the hypotheses of Lemma 3 are satisfied. The case $p = \infty$ is deferred to Sect. 3.1.

Let $1 \leq p < \infty$ and define the barrier

$$F(u, s) = F_p(u, s) = -\sum_i \log z_i - \sigma \sum_i \log s_i - \sum_i \log \tau_i \text{ where} \tag{25}$$

$$z_i = s_i^{2/p} - \sum_{j=1}^{d} [\overbrace{(D^{(j)}u + D^{(j)}g)_i}^{y^{(j)}}]^2, \qquad \tau_i = R - \omega_i s_i \text{ and} \tag{26}$$

$$\sigma = \sigma(p) = \begin{cases} 2 & \text{if } 1 \leq p < 2 \text{ and} \\ 1 & \text{if } p \geq 2. \end{cases} \tag{27}$$

**Lemma 6** *The function $F(u, t)$ is an $m(\sigma + 2)$-self-concordant barrier for the set*

$$Q = \{(u, s) \; : \; s_i \geq \|\nabla(u + g)|_{K_i}\|_2^p, \; s_i \geq 0 \text{ and } \max_i \omega_i s_i \leq R\}, \quad (28)$$

*The problem of minimizing $J(u)$ over $u \in V_h$ subject to the additional constraint that $\max_i \omega_i \|\nabla(u + g)|_{K_i}\|_2 \leq R$ is equivalent to:*

$$\min c^T x \text{ subject to } x \in Q \text{ where } c = \begin{bmatrix} -f \\ \omega \end{bmatrix}. \quad (29)$$

*Here, we have abused the notation and used the symbol $f$ for the vector whose $i$th component is $\int_\Omega f(x)\phi_i(x)\,dx$.*

**Proof** The functions $B_p(x, s) = -\log(s^{2/p} - x^T x)$ are $\sigma + 1$-self-concordant so $-\sum_i \log \tau_i + \sum_{i=1}^m B_p([\sum_k D_{i,k}^{(j)}(u + g)_k]_{j=1}^d, s_i)$ is $m(\sigma + 2)$ self-concordant, see Nesterov and Nemirovskii [28]. The rest is proved by inspection. □

From Lemma 3, it is tempting to use a bound such as $\|u\|_{X^p} < R$, i.e. $\sum_i \omega_u s_i \leq R$, but this leads to a dense Hessian $F_{ss}$. Instead, we have used the "uniform" bound:

$$\omega_i s_i \leq \sum_i \omega_i s_i = \int_\Omega s \leq R.$$

With this "looser" bound, the Hessian $F_{ss}$ is sparse.[3] Furthermore, by using the $R$ value from the a-priori estimate Lemma 3, one can ensure that $Q$ is non-empty and contains minimizing sequences for $J(u)$. Thus, put:

$$R \geq R^* = 2(1 + \|g\|_{X^p}^p) = 2 + 2\left(\sum_{j=1}^d [(D^{(j)}g)_i]^2\right)^{\frac{p}{2}}. \quad (30)$$

Set

$$\hat{s}_i = 1 + \left(\sum_{j=1}^d [(D^{(j)}g)_i]^2\right)^{\frac{p}{2}}; \text{ hence } \hat{s}_i \leq \frac{R}{2}. \quad (31)$$

In this way, $(0, \hat{s}) \in Q$.

**Lemma 7** *For all $(u, s) \in Q$,*

$$\tau_i \leq R, \quad s_i \leq \frac{R}{\omega_i}, \quad z_i \leq \left(\frac{R}{\omega_i}\right)^{\frac{2}{p}}. \quad (32)$$

---

[3] The alternative barrier term $-\log(R - \sum_i \omega_i s_i)$ corresponding to $\int_\Omega s < R$ leads to a rank-1 dense Hessian matrix so one could also use this more "natural" barrier by combining it with a Woodbury identity to invert the Hessian. This trades a slightly more complex implementation in exchange for a slightly "smaller" set $Q$.

**Proof** From $w^T s \geq 0$ and (26), we find $\tau \leq R$. From (26), we find $z_i \leq s_i^{2/p}$ and from $0 \leq \tau_i = R - \omega_i s$, we find $\omega_i s_i \leq R$.

We further find that:

$$\hat{\tau}_i \geq R/2, \quad \hat{s}_i \geq 1, \quad \hat{z}_i \geq 1. \tag{33}$$

The gradient of $F$ is:

$$F' = \begin{bmatrix} F_u \\ F_s \end{bmatrix} = \begin{bmatrix} \sum_j 2[D^{(j)}]^T \frac{y^{(j)}}{z} \\ -\frac{2}{p} \frac{1}{z} s^{2/p-1} - \frac{\sigma}{s} + \frac{\omega}{\tau} \end{bmatrix}, \tag{34}$$

where vector algebra is defined entrywise.

The Hessian of $F$ is

$$F'' = \begin{bmatrix} F_{uu} & F_{us} \\ F_{su} & F_{ss} \end{bmatrix} = \begin{bmatrix} F_{uu} & F_{us} \\ F_{us}^T & F_{ss} \end{bmatrix} \text{ where} \tag{35}$$

$$F_{uu} = 2\sum_{j=1}^{d} [D^{(j)}]^T Z^{-1} D^{(j)} + 4 \sum_{j,r=1}^{d} (Y^{(j)} D^{(j)})^T Z^{-2} (Y^{(r)} D^{(r)}), \tag{36}$$

$$F_{us} = -\frac{4}{p} \sum_{j=1}^{d} (Y^{(j)} D^{(j)})^T Z^{-2} S^{2/p-1}, \tag{37}$$

$$F_{ss} = -\frac{2}{p}\left(\frac{2}{p}-1\right) Z^{-1} S^{2/p-2} + \frac{4}{p^2} Z^{-2} S^{4/p-2} + \sigma S^{-2} + W^2 Z^{-2}, \tag{38}$$

$$S = \text{diag}(s), \ W = \text{diag}(\omega), \ Y = \text{diag}(y), \ Z = \text{diag}(z). \tag{39}$$

**Lemma 8** *Let $d_{\min}^2$ be the smallest eigenvalue of $\sum_{k=1}^{d} [D^{(k)}]^T D^{(k)}$ and assume $0 < h < 1$. Let $\omega_{\min} = \min_i \omega_i$, and similarly for $(z_F^*)_{\max}$, etc... The smallest eigenvalue $\lambda_{\min}$ of $F''(u_F^*, s_F^*)$ is bounded below by*

$$\lambda_{\min} \geq \min\left\{2(z_F^*)_{\max}^{-1} d_{\min}^2, \omega_{\min}^2 (z_F^*)_{\max}^{-2}\right\}. \tag{40}$$

**Proof** We consider the "Rayleigh quotient" $x^T F'' x / x^T x$, the extremal values of which are the extremal eigenvalues of $F''$. We put $x = \begin{bmatrix} v \\ w \end{bmatrix}$ so that

$$x^T F'' x = v^T F_{uu} v + 2v^T F_{us} w + w^T F_{ss} w$$

We use the Cauchy-Schwarz inequality together with Young's inequality to find that

$$2|v^T F_{us} w| = \frac{4}{p} \left| \left( \sum_{j=1}^{d} v^T (Y^{(j)} D^{(j)})^T Z^{-1} \right) \left( Z^{-1} S^{2/p-1} w \right) \right|$$

$$\leq \frac{8}{p} \sqrt{\sum_{j,r=1}^{d} v^T (Y^{(j)} D^{(j)})^T Z^{-2} (Y^{(r)} D^{(r)}) v} \sqrt{w^T S^{2/p-1} Z^{-2} S^{2/p-1} w}$$

$$\leq 4 \sum_{j=1,r}^{d} v^T (Y^{(j)} D^{(j)})^T Z^{-2} (Y^{(r)} D^{(r)}) v + \frac{4}{p^2} w^T S^{2/p-1} Z^{-2} S^{2/p-1} w.$$

Hence we find:

$$x^T F'' x \geq 2 \sum_{j=1}^{d} v^T [D^{(j)}]^T Z^{-1} D^{(j)} v + \left(\frac{2}{p} - 1\right) w^T \left(-\frac{2}{p} Z^{-1} S^{2/p-2}\right) w$$
$$+ w^T \sigma S^{-2} w + w^T W^2 Z^{-2} w.$$

We use that $F_s = 0$, which implies that $-\frac{2}{p} Z^{-1} S^{2/p-1} = T^{-1} W - \sigma S^{-1}$, where $T = \mathrm{diag}(\tau)$ and hence

$$x^T F'' x \geq 2 \sum_{j=1}^{d} v^T [D^{(j)}]^T Z^{-1} D^{(j)} v + \left(\frac{2}{p} - 1\right) w^T T^{-1} W S^{-1} w + w^T W^2 Z^{-2} w$$

$$\geq \|x\|_2^2 \min \left\{ 2 z_{\max}^{-1} d_{\min}^2, \omega_{\min}^2 z_{\max}^{-2} \right\}.$$

$\square$

A domain $\Omega$ is said to be of width $L$ when $\Omega \subset S$, where $S$ is a strip of width $L$. The Friedrichs inequality states that, for domains of width $L > 0$ and for $u \in W_0^{1,p}(\Omega)$, $\|u\|_{L^2} \leq L |u|_{H^1(\Omega)}$.

**Lemma 9** *Let $\Omega$ be a polytope of width $L < \infty$, and assume that the triangulation $T_h$, which depends on the grid parameter $0 < h < 1$, is quasi-uniform. Then, there is a constant $c_\Omega > 0$, which depends on $\Omega$ and the quasi-uniformity parameter $\rho$ of $T_h$, such that the smallest eigenvalue $d_{\min}^2$ of $\sum_k [D^{(k)}]^T D^{(k)}$ satisfies*

$$d_{\min}^2 \geq c_\Omega > 0. \tag{41}$$

**Proof** Consider the matrix $A = \sum_{k=1}^{d} [D^{(k)}]^T W D^{(k)}$ and note that

$$u^* A u \leq \omega_{\max} u^* \left(\sum_k [D^{(k)}]^T D^{(k)}\right) u \leq \frac{(\rho h)^d}{d!} u^* \left(\sum_k [D^{(k)}]^T D^{(k)}\right) u.$$

Furthermore, $u^T A u = |u|_{H^1}^2$, and the Friedrichs inequality states that $\|u\|_{L^2} \leq L |u|_{H^1}$. Furthermore, according to [32, Proposition 6.3.1], there is a constant $K_\Omega$ such that $u^T u \leq K_\Omega h^{-d} \|u\|_{L^2}^2$. Finally, we use that the quadrature weights $\{\omega_i\}$ are $\Theta(h^d)$ to find that $u^T u \leq C_\Omega h^{-d} |u|_{H^1}^2 = C_\Omega h^{-d} u^T A u \leq C_\Omega h^{-d} \frac{(\rho h)^d}{d!} u^* \left(\sum_k [D^{(k)}]^T D^{(k)}\right) u$, as required. $\square$

**Lemma 10** *Assume $T_h$ is quasi-uniform and that $R \geq 1$ and $0 < h \leq 1$. There is a constant $c'_\Omega > 0$, which depends on $\Omega$, such that*

$$\lambda_{\min} \geq c'_\Omega R^{-4} h^{6d}. \tag{42}$$

**Proof** Using (32) and (40), we arrive at:

$$\lambda_{\min} \geq \min \left\{ 2 \left( \frac{R}{\omega_{\min}} \right)^{-\frac{2}{p}} d_{\min}^2, \omega_{\min}^2 \left( \frac{R}{\omega_{\min}} \right)^{-\frac{4}{p}} \right\}.$$

Note that $R \geq R^* \geq 1 = |\Omega|^{-1} |\Omega| \geq |\Omega|^{-1} \omega_{\min}$ and hence $R/\omega_{\min} \geq |\Omega|^{-1}$ and

$$\lambda_{\min} \geq \min\{2c_\Omega, 1\} \max\{|\Omega|, 1\}^2 R^{-\frac{4}{p}} \omega_{\min}^{2+\frac{4}{p}}.$$

Since $R \geq 1$ and $\omega_{\min} \leq 1$ (because $h \leq 1$), we can find a lower bound by putting $p = 1$ in the exponents. Under the quasi-uniform hypothesis, all the quadrature weights are bounded below by $\omega_i \geq h^d/(d!)$, which yields (42). $\qquad\square$

**Lemma 11** *For $0 < h < 1$, assume $T_h$ is a quasi-uniform triangulation of $\Omega$. The number $n$ of vertices of $T_h$ inside $\Omega$ and the number $m$ of simplices in $T_h$ satisfy*

$$\frac{n}{d+1} \leq m \leq |\Omega| h^{-d} d!, \tag{43}$$

*where $|\Omega|$ is the volume of $\Omega$.*

**Proof** The inequality $n \leq (d+1)m$ follows from the fact that each of the $m$ simplices has precisely $d+1$ vertices; we may indeed have $n < (d+1)m$ since some vertices may be shared between multiple simplices. The upper bound for $m$ follows from $|\Omega| = \sum_{i=1}^m \omega_i \geq m\omega_{\min} \geq mh^d/(d!)$. $\qquad\square$

**Lemma 12** *Consider the point $\hat{x} = (0, \hat{s})$.*

$$\|\hat{x}\|_2 \leq C_\Omega^* h^{-1.5d} R, \text{ and } \|F'(\hat{x})\|_2 \leq C_\Omega^* h^{-1-1.5d} R(1 + \|g\|_{X^p}), \tag{44}$$

*where $C_\Omega^* < \infty$ is a constant that depends on $\Omega$ and the quasi-unifomity parameter $\rho$ of $T_h$.*

**Proof** From (31) we have

$$\|\hat{x}\|_2^2 = \sum_{i=1}^m \hat{s}_i^2 \leq m \left( \frac{R}{2\omega_{\min}} \right)^2. \tag{45}$$

Using (43), we obtain (44).

We now estimate $F'(\hat{x})$. Using (34), we find

$$\|F'(\hat{x})\|_2 \leq \|F_u\|_2 + \|F_s\|_2 \tag{46}$$

$$\leq \sum_j \|[D^{(j)}]^T \hat{Z}^{-1} D^{(j)} g\|_2 + \frac{2}{p} \hat{z}_{\min}^{-1} \|\hat{s}^{2/p-1}\|_2 + \sigma \|\hat{s}^{-1}\|_2 + \frac{1}{\hat{\tau}_{\min}} \|\omega\|_2. \tag{47}$$

We bound the first term as follows:

$$\sum_j \|[D^{(j)}]^T \hat{Z}^{-1} D^{(j)} g\|_2 \leq \left( \sum_j \|[D^{(j)}]^T \hat{Z}^{-1}\|_2^2 \right)^{1/2} \left( \sum_j \|D^{(j)} g\|_2^2 \right)^{1/2} \tag{48}$$

$$\leq \|\delta\|_2 \overbrace{\hat{z}_{\min}^{-1}}^{\leq 1} \omega_{\min}^{-1} \left( \sum_j \|D^{(j)} g\|_2^2 \right)^{1/2} \tag{49}$$

Here, $\delta_j^2 = \|D^{(j)}\|_2^2 \leq \omega_{\min}^{-1} \rho([D^{(j)}]^T W D^{(j)})$, where $\rho(\cdot)$ is the spectral radius. We estimate the spectral radius as follows:

$$w^T [D^{(j)}]^T W D^{(j)} w = \int_\Omega w_{x_j}^2 \, dx \leq |w|_{H^1}^2 \leq C_{IS}^2 h^{-2} \|w\|_{L^2}^2 \leq C_{IS}^2 [K_\Omega']^2 h^{2d-2} \|w\|_2^2,$$

where we have used the inverse Sobolev inequality $|w|_{H^1} \leq C_{IS} h^{-1} \|w\|_{L^2}$ for $w \in V_h$ (see e.g. Toselli and Widlund [37, Lemma B.27]) and the norm equivalence $\|w\|_{L^2} \leq K_\Omega' h^d \|w\|_2^2$. Thus,

$$\|\delta\|_2 \leq \omega_{\min}^{-1/2} \sqrt{d} C_{IS} [K_\Omega'] h^{d-1}.$$

Furthermore, using equivalence of $p-$norms in $m$ dimensions,

$$\left( \sum_j \|D^{(j)} g\|_2^2 \right)^{1/2} \leq \omega_{\min}^{-1/2} \left( \sum_{k=1}^m \omega_k \left[ \left( \sum_{j=1}^d (D^{(j)} g)_k^2 \right)^{1/2} \right]^2 \right)^{1/2} \tag{50}$$

$$\leq \omega_{\min}^{-1/2} m^{1/2} \left( \sum_{k=1}^m \omega_k \left[ \left( \sum_{j=1}^d (D^{(j)} g)_k^2 \right)^{1/2} \right]^p \right)^{1/p} \tag{51}$$

$$= \omega_{\min}^{-1/2} m^{1/2} \|g\|_{X^p}. \tag{52}$$

As a result,

$$\sum_j \|[D^{(j)}]^T \hat{Z}^{-1} D^{(j)} g\|_2 \leq \left( \omega_{\min}^{-1/2} \sqrt{d} C_{IS} [K_\Omega'] h^{d-1} \right) \left( \omega_{\min}^{-1/2} m^{1/2} \|g\|_{X^p} \right) \tag{53}$$

$$\leq C_{IS} K_\Omega'' h^{-1-0.5d} \|g\|_{X^p}. \tag{54}$$

From (31) and (33), we further estimate

$$\|\hat{s}^{2/p-1}\|_2 \leq \begin{cases} \sqrt{m}\hat{s}_{\max}^{2/p-1} \leq \sqrt{m}\left(\dfrac{R}{2\omega_{\min}}\right)^{2/p-1} & \text{if } 1 \leq p < 2 \\ \sqrt{m} & \text{if } p \geq 2 \end{cases} \tag{55}$$

$$\leq \sqrt{m}\left(\frac{R}{2\omega_{\min}}\right) \text{ and} \tag{56}$$

$$\|\hat{s}^{-1}\|_2 \leq \sqrt{m} \text{ and } \|\omega\|_2 \leq \sqrt{m}(\rho h)^d/(d!). \tag{57}$$

Hence,

$$\|F'(\hat{x})\|_2 \leq C_{IS}K''_{\Omega}h^{-1-0.5d}\|g\|_{X^p} \\ + \sqrt{m}\left(\frac{R}{2\omega_{\min}}\right) + 2\sqrt{m} + \frac{2}{R}\sqrt{m}(\rho h)^d/(d!).$$

$$\square$$

**Proof of Theorem 1** Using (42), we find

$$\|\hat{x}\|_{x_F^*}^* \leq \lambda_{\min}^{-1}\|\hat{x}\|_2 \leq \left(c'_{\Omega}R^{-4}h^{6d}\right)^{-1}\left(C^*_{\Omega}h^{-1.5d}R\right) = [c'_{\Omega}]^{-1}C^*_{\Omega}R^5h^{-7.5d}. \tag{58}$$

Also

$$\|F'(\hat{x})\|_{x_F^*}^* \leq \lambda_{\min}^{-1}\|F'(\hat{x})\|_2 \leq \left(c'_{\Omega}R^{-4}h^{6d}\right)^{-1}\left(C^*_{\Omega}h^{-1-1.5d}R(1+\|g\|_{X^p})\right) \tag{59}$$

$$= [c'_{\Omega}]^{-1}C^*_{\Omega}h^{-1-7.5d}R^5(1+\|g\|_{X^p}). \tag{60}$$

We substitute these estimates into (21) to get

$$N^* \leq 7.2\sqrt{4m}\left[2\log(4m)\right. \tag{61}$$

$$+ \log\left([c'_{\Omega}]^{-1}C^*_{\Omega}h^{-1-7.5d}R^5(1+\|g\|_{X^p})\right) \tag{62}$$

$$+ \left.\log\left([c'_{\Omega}]^{-1}C^*_{\Omega}R^5h^{-7.5d}\right) + \log(1/\epsilon)\right] \tag{63}$$

Using $m \leq |\Omega|d!h^{-d}$, we get

$$\leq 14.4\sqrt{|\Omega|d!h^{-d}}\left[\log\left(h^{-1-17d}R^5(1+\|g\|_{X^p})\epsilon^{-1}\right) + K^*\right]. \tag{64}$$

The estimates $R_{p=1}$, $R_{1<p<\infty}$ were obtained by starting from the estimates of Lemma 3, adding 1, and doubling them. Substituting these into $N^*$ produces $N_{p=1}$ and $N_{1<p<\infty}$. $\square$

## 3.1 The case $p = \infty$

Recall the $\infty$-Laplacian of (7). As in the $p = 1$ case, $J(u)$ is non-differentiable and may be unbounded below when $f$ is large. As per Lemma 3, assume that $L\|f\|_{L^1} < 1$ and set $R_{p=\infty} = \max_i \omega_i \left(2 + \frac{2\|g\|_{X^\infty(\Omega)}}{1 - L\|f\|_{L^1}}\right)$ and impose $\omega_i s \leq R_{p=\infty}$. The problem of minimizing $J(u)$ over $u \in V_h$ is equivalent to

$$\min s \text{ over } Q := \left\{(u, s) \ : \ s \geq \|\nabla(u + g)|_{K_i}\|_2 \ \forall i, \text{ and } R \geq \omega_i s\right\} \quad (65)$$

We notice that this definition of $Q$ coincides with the definition (28) with $p = 1$ subject to the additional restriction that $s_1 = \ldots = s_m$ and subsequently dropping the index $i$ from $s_i$. As a result, one can obtain a barrier for $Q$ by taking the barrier (25) with $p = 1$ on the subspace of constant valued $s$ vectors, hence the barrier $F_\infty$ and its derivatives are

$$F_\infty(u, s) = F_1(u, se), \ F'_\infty(u, s) = \mathcal{E}F'_1(u, se), \ F''_\infty(u, s) = \mathcal{E}F''_1(u, se)\mathcal{E}^T, \quad (66)$$

$$\text{where } e = \begin{bmatrix} 1 \\ \vdots \\ 1 \end{bmatrix}, \ \mathcal{E} = \begin{bmatrix} I & O \\ O & e^T \end{bmatrix}. \quad (67)$$

The starting point for the optimization is $(\hat{u}, \hat{s})$ with $\hat{u} = 0$ and

$$\hat{s} = 1 + \max_i \left(\sum_{j=1}^d [(D^{(j)}g)_i]^2\right)^{\frac{1}{2}}. \quad (68)$$

**Theorem 3** *With the notation as in Theorem 1, let*

$$R = R_{p=\infty} = \max_i \omega_i \left(2 + \frac{2\|g\|_{X^\infty(\Omega)}}{1 - L\|f\|_{L^1}}\right) \quad (69)$$

*and assume $p = \infty$, $L\|f\|_{L^1} < 1$. The barrier method to solve (65) requires at most $N_{p=\infty}$ Newton iterations, where*

$$N_{p=\infty} \leq 14.4\sqrt{|\Omega|h^{-d}d!}\left[\log\left(h^{-1-6.5d}\left(2 + \frac{2\|g\|_{X^\infty(\Omega)}}{1 - L\|f\|_{L^1}}\right)^5 \epsilon^{-1}\right) + K^*\right]. \quad (70)$$

*The computational complexity as a function of the number $n$ of grid points (and freezing all other parameters) is $O(\sqrt{n}\log n)$.*

**Proof** The proof of Theorem 3 follows the same logic as that of Theorem 1, so we merely sketch it here. First, (34) and (35) are replaced by:

$$F' = \begin{bmatrix} F_u \\ F_s \end{bmatrix} = \begin{bmatrix} \sum_j 2[D^{(j)}]^T \frac{y^{(j)}}{z} \\ -2s \sum_j \frac{1}{z_j} - \frac{m\sigma}{s} + \sum_j \frac{\omega_j}{\tau_j} \end{bmatrix}, \qquad (71)$$

$$F'' = \begin{bmatrix} F_{uu} & F_{us} \\ F_{su} & F_{ss} \end{bmatrix} = \begin{bmatrix} F_{uu} & F_{us} \\ F_{us}^T & F_{ss} \end{bmatrix} \text{ where} \qquad (72)$$

$$F_{uu} = 2\sum_{j=1}^{d} [D^{(j)}]^T Z^{-1} D^{(j)} + 4\sum_{j,r=1}^{d} (Y^{(j)} D^{(j)})^T Z^{-2} (Y^{(r)} D^{(r)}), \qquad (73)$$

$$F_{us} = -4\sum_{j=1}^{d} (Y^{(j)} D^{(j)})^T z^{-2} s, \qquad (74)$$

$$F_{ss} = -2\sum_j z_j^{-1} + 4\sum_j z_j^{-2} s^2 + \sigma m s^{-2} + \sum_j \omega_j^2 z_j^{-2}. \qquad (75)$$

The proof of (40) holds (changing what must be changed), ending with

$$x^T F'' x \geq \|x\|_2^2 \min\{2 z_{\max}^{-1} d_{\min}^2, \sum_k \omega_k^2 z_k^{-2}\}, \qquad (76)$$

which is slightly stronger than (40).

The estimate (44) also holds verbatim. The estimate for $\|\hat{x}\|_2$ is by inspection of (68) and (69). The estimate for $\|F_u(\hat{x})\|_2$ is identical to the proof of (44), and $|F_s(\hat{x})|$ is estimated as follows:

$$|F_s| \leq 2\left(\frac{R}{2\omega_{\min}}\right) z_{\min}^{-1} m + \sigma m + |\Omega| \tau_{\min}^{-1} \leq C_\Omega R h^{-2d}, \qquad (77)$$

where we have used (16), $1 \leq \hat{s} \leq R/(2\min\omega_i)$, $\hat{\tau}_i \geq R/(2\min\omega_i)$, and $\hat{z}_j \geq 1$, and $C_\Omega$ is some constant that depends only on $\Omega$. Thus,

$$\|\hat{x}\|_{x_F^*}^* \leq \lambda_{\min}^{-1} \|\hat{x}\|_2 \leq [c_\Omega']^{-1} C_\Omega^* R^5 h^{-7.5d},$$

see (58). Futhermore,

$$\|F'(\hat{x})\|_{x_F^*}^* \leq \lambda_{\min}^{-1} \|F'(\hat{x})\|_2 \qquad (78)$$

$$\leq \left(c_\Omega' R^{-4} h^{6d}\right)^{-1} \left(C_{IS} K_\Omega'' h^{-1-0.5d} \|g\|_{X^\infty} + C_\Omega R h^{-2d}\right) \qquad (79)$$

$$\leq K_\Omega''' h^{-1-3d} \left(\frac{\|g\|_{X^\infty}}{1 + \|f\|_{X^1}}\right), \qquad (80)$$

yielding the final estimate

$$N^* \leq 7.2\sqrt{4m} \left[ 2\log(4m) \right. \tag{81}$$

$$+ \log\left( K_\Omega''' h^{-1-3d} \left( \frac{\|g\|_{X^\infty}}{1 + \|f\|_{X^1}} \right) \right) \tag{82}$$

$$+ \log\left( [c_\Omega']^{-1} C_\Omega^* R^5 h^{-7.5d} \right) \tag{83}$$

$$+ \left. \log(\epsilon^{-1}) \right], \tag{84}$$

as required. □

## 3.2 Implementation notes

In principle, the vector $(f_i)$ is defined by $f_i = \int_\Omega f \phi_i$; we have not analyzed an inexact scheme for computing these integrals. If $f$ is assumed to be a suitable finite element space (e.g. piecewise constant), then these integrals can be computed exactly from the same quadrature we use on the diffusion term. In addition, we can then compute $\|f\|_{L^q}$ exactly by quadrature since $|f|^q$ is also piecewise constant. Assuming $g$ is piecewise linear, the quantities $|\Omega|$ and $\|g\|_{X^p}^p$ can be computed exactly, see (30). Thus, one can compute $R_{1<p<\infty}$, $R_{p=1}$, etc... exactly.

In the strong form (1), the function $g$ is given on $\partial\Omega$ (i.e. it is a trace), but in the variational form (2), the function $g$ has domain $\Omega$. Regarding $v = u + g$ as the solution, the choice of $g|_\Omega$ doe not affect the value of $v$, provided that $g|_{\partial\Omega}$ is fixed. The simplest way to choose $g|_\Omega$ as a piecewise linear function on $T_h$ is to set all nodal values to 0 inside of $\Omega$, but this is a somewhat "rough" prolongation that is furthermore dependent on $h$. Using such a prolongation of $g$ causes the estimates (5) and (6) to become dependent on $h$ where $g$ appears. In order to avoid this dependence on $h$, one can proceed in one of two ways. First, if the meshes $T_h$ are all included in one coarse mesh $T_{h_0}$, then one can do the prolongation on $T_{h_0}$ and use the same prolongation on all $T_h$.

Another method is to solve the linear Laplacian with boundary conditions $g|_{\partial\Omega}$ on the mesh $T_h$. This choice of $g|_\Omega$ does vary slightly with $h$ but it converges to the continuous solution as $h \to 0$. Furthermore, this choice of $g$ minimizes $\|g\|_{X^2} = |g|_{H^1}$ so it may result in a smaller value of $R$ than prolongation by truncation. We call this choice of $g$ the **discrete harmonic prolongation** of $g|_{\partial\Omega}$ to $\Omega$. We use the discrete harmonic prolongation in all our numerical experiments.

## 4 Numerical experiments

We consider the *p*-Laplacian for $p = 1, 1.1, 1.2, 1.5, 2, 3, 4, 5, \infty$ for a square domain subject to Dirichlet boundary conditions and where the forcing $f = 0$, see Fig. 1. For the boundary conditions $g$, we have taken the piecewise linear interpolant of the trace $(1_X(x, y))|_{\partial\Omega}$ of the characteristic function $1_X(x, y)$, where $X$ is the set

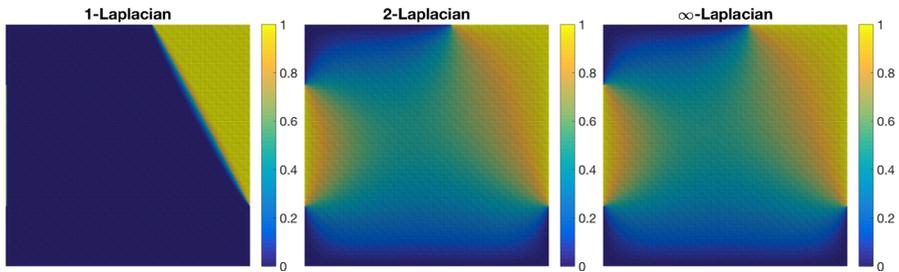

**Fig. 1** Solving the *p*-Laplacian for $p = 1, 2, \infty$ with the same boundary conditions $g$ and zero forcing $f = 0$ on a $200 \times 200$ grid. Because of the zero forcing, the minimum and maximum principles hold, which provides some protection against the near-discontinuities in the boundary data, e.g. when $p = \infty$

$X = (\{0\} \times [0.25, 0.75]) \cup ([0.6, 1] \times [0.25, 1])$, which we approximate on the discrete grid by piecewise linear elements. Note that this creates very challenging numerical and functional analytical problems, e.g. the trace of $W^{1,\infty}$ functions are also $W^{1,\infty}$ so the $\infty$−Laplacian here is solving a problem approximating one outside the usual trace space. The forcing $f = 0$ means that solutions must satisfy minimum and maximum principles, and so the solution is always between the extremal 0 and 1 boundary values for all values of $p$ and all $x \in \Omega$. The zero forcing provides some "protection" against the "bad" boundary data.

We have varied the number $n$ of grid points from $n = 16$ (a $4 \times 4$ grid) up to $n = 40,000$ ($200 \times 200$) and in all cases, solved to a tolerance of $\epsilon \approx 10^{-6}$. We have reported the number of Newton iterations required for convergence in Table 1. This detailed table reveals those values of $\kappa$, $n$, $p$ that failed to converge within five minutes. Most of these convergence failures are due to purely numerical problems. Indeed, we have noted in the introduction that when $p$ is large, minimizing $J(u)$ is intrinsically challenging because it exhausts the accuracy of double precision floating point. Thus, the difficulty in solving $p$-Laplacians accurately for large $p$ is not particular to our algorithm but indeeds affects all algorithms for solving $p$-Laplacians. MATLAB has also issued warnings that the Hessian was singular to machine precision, for large values of $p$ and $n$.

The scaling properties of our algorithms are not immediately obvious from Table 1. In order to visualize the scaling properties of our algorithms, we have sketched the iterations counts of Table 1 in Fig. 2. Note that both axes are in logarithmic scale, so straight lines of slope $\alpha$ correspond to $O(n^{\alpha})$ scaling. We see that the short step algorithm of Sect. 2.3 requires the largest number of Newton iterations to converge (blue lines). This is consistent with experience in convex optimization. For this reason, we were not able to solve larger problems with the short-step algorithm. The scaling of the short-step algorithm is consistent with the theoretical prediction $O(\sqrt{n} \log n)$ of Theorems 1 and 3.

The long step algorithms (black lines) all require fewer Newton steps than the short step algorithm, even though the theoretical estimate $O(n \log n)$ for long step algorithms is worse than for short step algorithms. This is a well-known phenomenon, and in practice, long step algorithms perform better, as is the case here.

**Table 1** Newton iteration counts for various problem sizes $n$, various step strategies $\kappa$ and various values of $p$

| $\kappa$ | $n = 16$ | 36 | 64 | 100 | 169 | 400 | 900 | 2500 | 5625 | 10,000 | 22,500 | 40,000 |
|---|---|---|---|---|---|---|---|---|---|---|---|---|
| $p = 1.0$ | | | | | | | | | | | | |
| Short | 985 | 2165 | 3163 | 4392 | 6086 | 10, 512 | 17, 192 | — | — | — | — | — |
| $\kappa = 2$ | 128 | 141 | 164 | 171 | 190 | 211 | 243 | 289 | 357 | 392 | 496 | 533 |
| $\kappa = 3$ | 88 | 98 | 136 | 118 | 148 | 165 | 194 | 237 | 343 | 535 | 825 | — |
| $\kappa = 4$ | 60 | 68 | 114 | 96 | 124 | 138 | 160 | 340 | 496 | 720 | 877 | — |
| $\kappa = 7$ | 68 | 67 | 136 | 99 | 252 | 255 | 384 | 464 | 833 | 1113 | — | — |
| Adaptive | 70 | 76 | 181 | 106 | 163 | 169 | 216 | 270 | 275 | 315 | 405 | 449 |
| $p = 1.1$ | | | | | | | | | | | | |
| Short | 954 | 1844 | 2677 | 3674 | 5094 | 8803 | 14, 510 | — | — | — | — | — |
| $\kappa = 2$ | 128 | 141 | 156 | 166 | 170 | 195 | 225 | 247 | 272 | 279 | 312 | 327 |
| $\kappa = 3$ | 89 | 96 | 103 | 113 | 116 | 130 | 157 | 185 | 211 | 229 | 244 | 274 |
| $\kappa = 4$ | 61 | 77 | 84 | 97 | 102 | 115 | 145 | 164 | 182 | 195 | 227 | 256 |
| $\kappa = 7$ | 68 | 69 | 74 | 91 | 88 | 104 | 143 | 172 | 190 | 237 | 314 | 391 |
| Adaptive | 60 | 81 | 68 | 97 | 80 | 103 | 190 | 186 | 243 | 210 | 231 | 274 |
| $p = 1.2$ | | | | | | | | | | | | |
| Short | 938 | 1779 | 2626 | 3566 | 4991 | 8589 | 14, 099 | — | — | — | — | — |
| $\kappa = 2$ | 128 | 141 | 156 | 162 | 169 | 193 | 218 | 244 | 259 | 265 | 291 | 295 |
| $\kappa = 3$ | 87 | 102 | 107 | 110 | 116 | 126 | 138 | 162 | 176 | 189 | 199 | 210 |
| $\kappa = 4$ | 60 | 67 | 74 | 86 | 89 | 98 | 114 | 124 | 139 | 152 | 168 | 176 |
| $\kappa = 7$ | 60 | 71 | 72 | 80 | 83 | 97 | 110 | 118 | 129 | 136 | 140 | 170 |
| Adaptive | 54 | 57 | 64 | 70 | 74 | 84 | 105 | 113 | 127 | 198 | 147 | 204 |

**Table 1** continued

| $\kappa$ | $n=16$ | 36 | 64 | 100 | 169 | 400 | 900 | 2500 | 5625 | 10,000 | 22,500 | 40,000 |
|---|---|---|---|---|---|---|---|---|---|---|---|---|
| $p=1.5$ | | | | | | | | | | | | |
| Short | 922 | 1731 | 2579 | 3484 | 4896 | 8405 | 13,754 | — | — | — | — | — |
| $\kappa=2$ | 129 | 146 | 154 | 163 | 174 | 193 | 220 | 249 | 266 | 273 | 287 | 294 |
| $\kappa=3$ | 86 | 102 | 106 | 111 | 116 | 132 | 142 | 161 | 175 | 186 | 194 | 201 |
| $\kappa=4$ | 58 | 71 | 73 | 84 | 88 | 102 | 116 | 126 | 131 | 148 | 156 | 162 |
| $\kappa=7$ | 59 | 64 | 77 | 80 | 83 | 98 | 108 | 117 | 120 | 128 | 126 | 139 |
| Adaptive | 52 | 58 | 66 | 69 | 72 | 85 | 96 | 111 | 110 | 120 | 123 | 130 |
| $p=2.0$ | | | | | | | | | | | | |
| Short | 829 | 1543 | 2305 | 3102 | 4356 | 7440 | 12,123 | 22,082 | — | — | — | — |
| $\kappa=2$ | 135 | 151 | 158 | 165 | 174 | 190 | 220 | 256 | 274 | 287 | 296 | 309 |
| $\kappa=3$ | 93 | 100 | 108 | 113 | 120 | 129 | 146 | 169 | 184 | 192 | 201 | 209 |
| $\kappa=4$ | 51 | 64 | 70 | 71 | 84 | 89 | 99 | 109 | 119 | 129 | 139 | 142 |
| $\kappa=7$ | 62 | 74 | 80 | 81 | 98 | 113 | 125 | 132 | 145 | 148 | 158 | 160 |
| Adaptive | 59 | 63 | 69 | 73 | 79 | 86 | 99 | 107 | 113 | 118 | 126 | 134 |
| $p=3.0$ | | | | | | | | | | | | |
| Short | 857 | 1624 | 2451 | 3329 | 4701 | 8127 | 13,372 | — | — | — | — | — |
| $\kappa=2$ | 142 | 161 | 172 | 184 | 195 | 208 | 249 | 299 | 315 | 331 | 348 | 363 |
| $\kappa=3$ | 99 | 113 | 122 | 129 | 136 | 144 | 164 | 193 | 212 | 223 | 239 | 249 |
| $\kappa=4$ | 56 | 65 | 74 | 79 | 87 | 102 | 110 | 128 | 136 | 157 | 163 | 163 |
| $\kappa=7$ | 67 | 75 | 84 | 90 | 102 | 116 | 127 | 138 | 154 | 156 | 165 | 175 |
| Adaptive | 84 | 70 | 76 | 112 | 88 | 94 | 109 | 118 | 129 | 135 | 145 | 147 |

Table 1 continued

| κ | n = 16 | 36 | 64 | 100 | 169 | 400 | 900 | 2500 | 5625 | 10,000 | 22,500 | 40,000 |
|---|---|---|---|---|---|---|---|---|---|---|---|---|
| **p = 4.0** | | | | | | | | | | | | |
| Short | 907 | 1765 | 2697 | 3696 | 5260 | 9197 | 15,265 | — | — | — | — | — |
| κ = 2 | 153 | 179 | 192 | 205 | 217 | 236 | 273 | 346 | 373 | 413 | 433 | — |
| κ = 3 | 105 | 122 | 132 | 144 | 153 | 167 | 189 | 224 | 243 | 289 | — | — |
| κ = 4 | 57 | 75 | 82 | 84 | 97 | 115 | 127 | 142 | 158 | 219 | 186 | 212 |
| κ = 7 | 69 | 84 | 94 | 99 | 108 | 128 | 141 | 159 | 174 | 178 | 193 | — |
| Adaptive | 66 | 77 | 80 | 91 | 100 | 110 | 125 | 137 | 148 | 157 | 512 | — |
| **p = 5.0** | | | | | | | | | | | | |
| Short | 970 | 1927 | 2974 | 4102 | 5872 | 10,352 | 17,280 | — | — | — | — | — |
| κ = 2 | 162 | 194 | 213 | 226 | 245 | 268 | 303 | 414 | — | — | — | — |
| κ = 3 | 112 | 132 | 145 | 159 | 169 | 188 | 214 | — | — | — | — | — |
| κ = 4 | 63 | 78 | 90 | 95 | 111 | 132 | 143 | 176 | — | — | — | — |
| κ = 7 | 74 | 91 | 102 | 109 | 119 | 146 | 164 | — | — | — | — | — |
| Adaptive | 67 | 82 | 91 | 100 | 109 | 123 | 136 | 430 | — | — | — | — |
| **p = ∞** | | | | | | | | | | | | |
| Short | 560 | 841 | 1191 | 1560 | 2145 | 3614 | 5884 | 10,809 | — | — | — | — |
| κ = 2 | 96 | 103 | 111 | 120 | 130 | 137 | 160 | 771 | 309 | 489 | — | — |
| κ = 3 | 69 | 73 | 78 | 86 | 91 | 95 | 113 | 451 | 536 | — | — | — |
| κ = 4 | 55 | 53 | 52 | 64 | 58 | 97 | 78 | 2606 | 3489 | — | — | — |
| κ = 7 | 54 | 48 | 57 | 62 | 63 | 75 | 1272 | 6989 | — | — | — | — |
| Adaptive | 50 | 62 | 66 | 62 | 73 | 83 | 120 | 150 | 276 | 191 | 302 | 341 |

The symbol — indicates that the algorithm failed to converge by exceeding the 5 minutes time limit

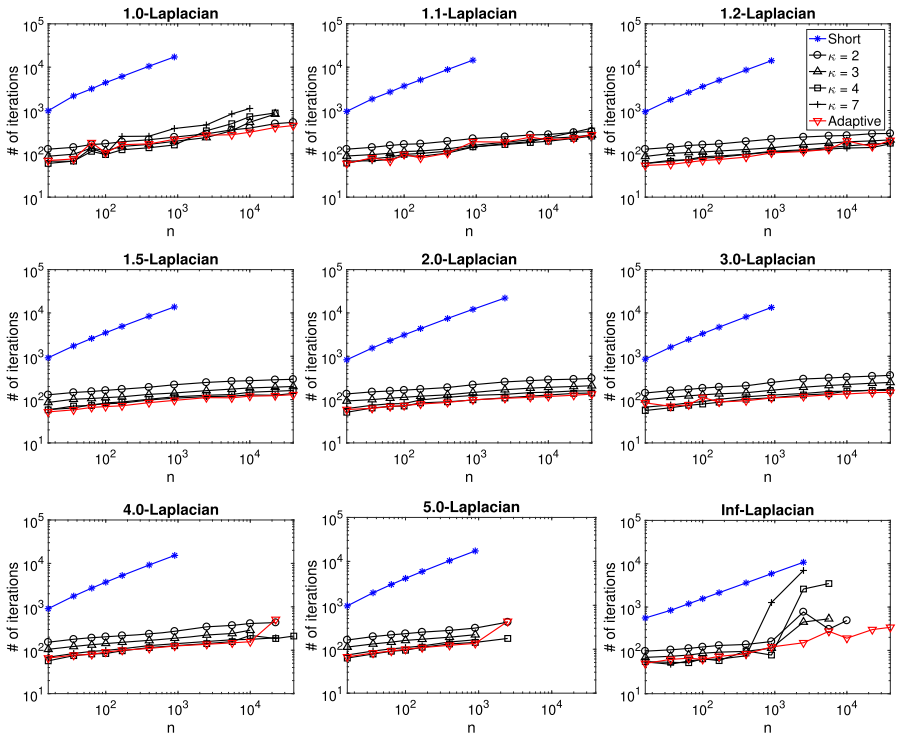

**Fig. 2** The number of Newton iterations for various grid sizes $n$ and parameters $p$ and step sizes $\kappa$

In Fig. 2, most of the black curves are approximately straight lines, indicating $O(n^\alpha)$ scaling, but there are notable exceptions when $p = 1$ or $p = \infty$, especially when $\kappa$ is also large. By contrast, the adaptive step size algorithm (red lines), with $\kappa_0 = 10$, is seen to be the best algorithm in most cases, and these red lines are much straighter than the black lines. We denote by $N_p(n)$ the number of iterations required for a certain value of $p$ and problem size $n$ for the adaptive step size algorithm. We have fitted straight lines to the red curves of Fig. 2 in the least-squares sense and obtained the following approximations:

$$
\begin{array}{cccccccccc}
p = & 1.0 & 1.1 & 1.2 & 1.5 & 2.0 & 3.0 & 4.0 & 5.0 & \infty \\
N_p(n) \approx & 62n^{0.18} & 33n^{0.21} & 31n^{0.17} & 43n^{0.11} & 47n^{0.10} & 60n^{0.09} & 30n^{0.22} & 17n^{0.36} & 18n^{0.28}
\end{array}
$$

Thus, it seems like the adaptive scheme requires about $O(n^{\frac{1}{4}})$ Newton iterations, regardless of the value of $p$.

Note that the case $p = 2$ is a linear Laplacian that can be computed by solving a single linear problem. When we embed this linear problem into the machinery of convex optimization, the overall algorithm is very inefficient since it may require hundreds of linear solves. We are including this test case for completeness, not as a recommendation.

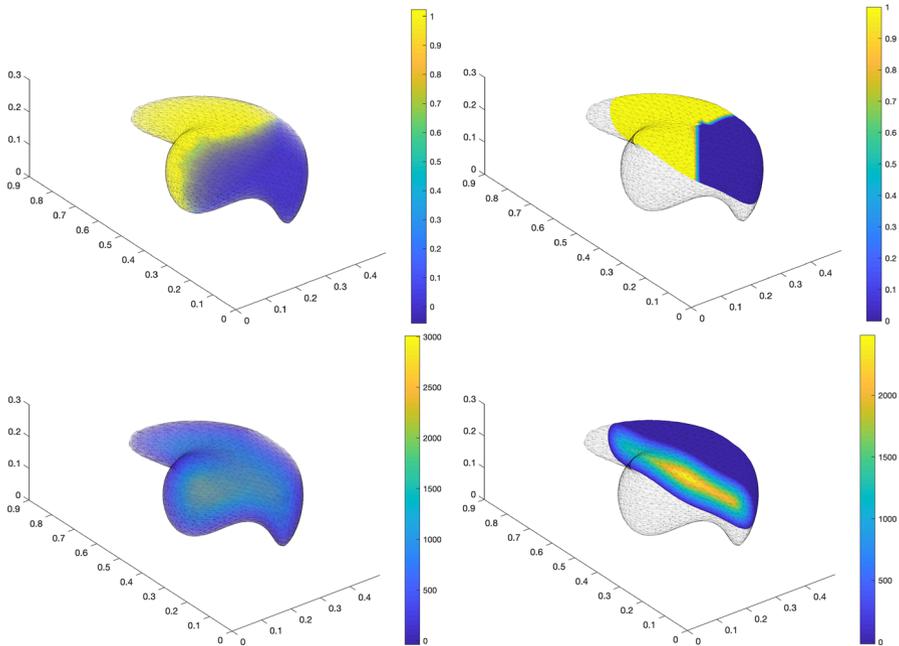

**Fig. 3** Solving the 1-Laplacian (top row) and $\infty$-Laplacian (bottom row) in 3d. The left column shows the solutions on the whole volumetric domain $\Omega$ with transparency, while the right column shows a slice through $\Omega$ of the same solutions with opaque colors

### 4.1 3d experiments

Consider the following function:

$$\phi = \frac{9}{20} - \sqrt{\left(x^2 + y^2\right)\left(1/10 + (|x - \cos(y)|)^2\right) + \left(z + \frac{3\,e^{-x}}{25}\right)^2}. \quad (85)$$

We define the "spaceship domain" $\tilde{\Omega} = \{(x, y, z) \in \mathbb{R}^3 \ : \ \phi > 0\}$; this domain is slightly rescaled so that it is aesthetically pleasing. In practice, the domain $\tilde{\Omega}$ is approximated on a discrete grid with a tetrahedral mesh $T_h$ and the corresponding polyhedral approximation $\Omega$ of $\tilde{\Omega}$. On this tetrahedral mesh, we solve the $p-$Laplacian with forcing $f = 1$ with $p \in \{1, \infty\}$. The boundary values $g$ are the indicating function of the set $\{y > 0.45\}$, as approximated by a piecewise linear function on the finite element grid. This problem features $n = 11,224$ unknowns and $m = 47,956$ elements. The solutions are displayed in Fig. 3.

For these problems, the solution of the 1-Laplacian seems to approximate the indicating function of $\{y > 0.45\}$, as expected. However, the solution of the $\infty$-Laplacian is very large (exceeding 2,000 somewhere in the middle of the spaceship). This is because the traces of $W^{1,\infty}(\Omega)$ functions are in $W^{1,\infty}(\partial\Omega)$ but our boundary data $g$ is a piecewise linear approximation of a discontinuous trace with jumps (an indicating function), an hence $\|g\|_{X^\infty}$ is very large and so is the solution $u + g$. The 1-Laplacian is better able to tolerate the boundary data $g$ with (near)-jumps because the trace of a $W^{1,1}(\Omega)$ function is merely $L^1(\partial\Omega)$, thus allowing jumps.

The solution for the $p = 1$-Laplacian seems very close to what one would obtain if one were to put $f = 0$ instead of $f = 1$. This is not surprising, because the 1-Laplacian is a linear program and the solutions of linear programs change in discrete steps when the forcing changes continuously. For example, the unique minimizer of $\tilde{J}(x) = |x| + fx$ ($x \in \mathbb{R}$) is $x = 0$ whenever $|f| < 1$, and switches to "undefined" (or $\pm\infty$) when $|f| > 1$ because then $\tilde{J}$ is unbounded below.

For the $p = \infty$-Laplacian, the solution $u + g$ is a large positive bump because $f > 0$ and there is a minimum principle stating that the minimum must of $u + g$ be on the boundary $\partial\Omega$. When one takes $f < 0$ instead, the solution $u + g$ is a large negative bump because in that scenario, $u + g$ satisfies a maximum principle. In the 2d experiments, the $\infty$-Laplacian did not develop large bumps because the boundary data was between 0 and 1 and the forcing was 0. This meant that $u + g$ had to satisfy both minimum and maximum principles, and $u$ was constrained by $0 \leq u + g \leq 1$, preventing the formation of large bumps in the solution.

# 5 Conclusions and outlook

We have presented new algorithms for solving the $p$-Laplacian efficiently for any given tolerance and for all $1 \leq p \leq \infty$. We have proven that our algorithms compute a solution to any given tolerance in polynomial time, using $O(\sqrt{n} \log n)$ Newton iterations, and an adaptive stepping variant converges in $O(\sqrt{n} \log^2 n)$ Newton iterations. We have confirmed these scalings with numerical experiments. We have further shown by numerical experiments that the adaptive step variant of the barrier method converges much faster than the short-step variant for the $p$-Laplacian and also usually faster than long-step barrier methods, thus achieving the practical speedup of long-step algorithms while avoiding the $O(n \log n)$ worst-case behavior of long-step algorithms. We have numerically estimated that the adaptive step algorithm requires $O(n^{\frac{1}{4}})$ Newton iterations across all values of $1 \leq p \leq \infty$. We have observed numerical difficulties for $p \geq 5$, which are expected since large powers exhaust the accuracy of double precision floating point arithmetic; this difficulty is not specific to our algorithm but is inherent to the $p$-Laplacian for large values of $p$. Our algorithms are particularly attractive when $p \approx 1$ and $p = \infty$, where there are no other algorithms that are efficient at all tolerances.